\documentclass[12pt]{article}
\usepackage[applemac]{inputenc}
\usepackage{amsmath,amssymb,amsfonts}
\usepackage{geometry}
\newtheorem{theorem}{Theorem}[section]
\newtheorem{lemma}[theorem]{Lemma}

\newtheorem{definition}[theorem]{Definition}

\newtheorem{problem}[theorem]{Problem}
\newtheorem{remark}[theorem]{Remark}

\newenvironment{myenumerate}{

\begin{enumerate}}{\end{enumerate}}
\newcommand{\dproof}{\noindent {Proof.} \quad}
\newcommand{\fproof}{\hfill $\square$ \bigskip}

\numberwithin{equation}{section}

\def\E{\mathbb{E}}

\def\1B{\text{1\!\!I}}

\def\r{\rangle}
\def\<{\langle}
\def\>{\rangle}

\def\R{\mathbb{R}}

\def\P{\mathbb{P}}
\def\Y{\mathcal{Y}}

\def\L{\mathcal{L}}

\def\X{\widehat{X}}

\begin{document}

\title{A time-fractional Kalman filter}
\author{Olfa Draouil$^{1}$, Rahma Yasmina Moulay Hachemi$^{2}$\\ 
        Bernt \O ksendal$^{3}$\\   
        Abderrahmen Aliane ${4}$}
\date{25 August 2026}
\maketitle

\footnotetext[1]{Department of Mathematics, University of Tunis El Manar, Tunisia\newline
Email: olfa.draouil@fst.utm.tn}

\footnotetext[2]{%
Department of Mathematics, University of Oslo, Norway. \newline
Email: yasmin.moulayhachemi@yahoo.com}

\footnotetext[3]{%
Department of Mathematics, University of Oslo, Norway. \\
Email: oksendal@math.uio.no.}

\footnotetext[4]{%
Department of Psychology, Educational sciences and Orthophony, Abou Bekr Belkaid University - Algeria. \\
Email: aliane.abderrahmen@gmail.com.}

\begin{abstract}
\small
We study a linear filtering problem in which the signal process is described
by a time-fractional linear stochastic differential equation driven by
Brownian motion. We derive a stochastic integral equation for the conditional
mean alongside a Riccati--Volterra type integral equation for the mean-square
error function. As a core application, we introduce a time-fractional
state-estimation framework for modelling learning trajectories in children
with developmental dyscalculia.
\end{abstract}

\vspace{0.2cm}



\section{Introduction}
The Kalman filter, introduced in the 1960s by R.E. Kalman \cite{K}, revolutionized the field of signal processing and control theory by providing an efficient recursive solution to the linear quadratic estimation problem. The Kalman filter consists of a series of mathematical equations that offer an effective recursive method for estimating the state of a process while minimizing the mean squared error. Traditional applications of the Kalman filter typically involve time-evolving signals and observations influenced by temporal Brownian motion. For more details about linear filtering, we refer to K\"orezli\v{o}glu \cite{Ko}, \O ksendal \cite{Oksendal13}, Wong \cite{Wo}, and for nonlinear filtering, we refer to Crisan \& Rozovskii \cite{CR}, Jazwinski \cite{J}, K\"orezli\v{o}glu et al. \cite{KMS}. 

Our proposed application is the estimation of cognitive learning
trajectories from repeated educational observations. Cognitive ability is not
observed directly, and present performance depends on a long history of
learning, forgetting, and previous difficulties. This motivates the use of a
fractional hidden-state model. In Section~\ref{sec:dyscalculia-example}, we
develop a fractional Kalman filtering framework for hidden cognitive-state
estimation in children with developmental dyscalculia.

\section{Preliminaries}
\label{sec:preliminaries}

This section recalls the fractional-calculus definitions and transform
identities used throughout the paper.

\subsection{Mittag--Leffler functions}

\begin{definition}[Two-parameter Mittag--Leffler function]
For $z,\alpha,\beta\in\mathbb C$ with
$\operatorname{Re}(\alpha)>0$ and $\operatorname{Re}(\beta)>0$, define
\begin{equation}
E_{\alpha,\beta}(z)
=\sum_{k=0}^{\infty}\frac{z^k}{\Gamma(\alpha k+\beta)},
\end{equation}
where $\Gamma$ denotes the Gamma function.
\end{definition}

\begin{definition}[One-parameter Mittag--Leffler function]
For $z,\alpha\in\mathbb C$ with $\operatorname{Re}(\alpha)>0$, define
\begin{equation}
E_\alpha(z)=E_{\alpha,1}(z)
=\sum_{k=0}^{\infty}\frac{z^k}{\Gamma(\alpha k+1)}.
\end{equation}
\end{definition}

\subsection{Fractional integrals and the Caputo derivative}

\begin{definition}[Riemann--Liouville fractional integral]
For $\alpha>0$, the fractional integral of a locally integrable function $f$
is
\begin{align}
(I^\alpha f)(t)
=\frac{1}{\Gamma(\alpha)}\int_0^t(t-s)^{\alpha-1}f(s)ds.
\end{align}
\end{definition}

\begin{definition}[Caputo fractional derivative]
Let $\alpha>0$ and let $n=\lceil\alpha\rceil$. The Caputo derivative is
defined by
\begin{align}\label{caputo1}
D_C^\alpha f(t)
&=I^{n-\alpha}f^{(n)}(t)\nonumber\\
&=\frac{1}{\Gamma(n-\alpha)}
  \int_0^t(t-s)^{n-\alpha-1}f^{(n)}(s)ds,
  \qquad n-1<\alpha<n,
\end{align}
while $D_C^n f=f^{(n)}$ when $\alpha=n$.
\end{definition}

For example, if $f(t)=t$ and $0<\alpha<1$, then
\begin{align}
D_C^\alpha f(t)
=\frac{t^{1-\alpha}}{\Gamma(2-\alpha)}.
\end{align}
In particular, $D_C^{1/2}t=2\sqrt{t}/\sqrt{\pi}$.

\subsection{Laplace-transform identities}

For a function $f$ of suitable growth, let
\begin{equation}
(Lf)(s)=\int_0^\infty e^{-st}f(t)dt=: \widetilde f(s).
\end{equation}
For $0<\alpha\leq1$, the identities used below are
\begin{align}
L[D_C^\alpha f](s)
&=s^\alpha(Lf)(s)-s^{\alpha-1}f(0),                         \label{L1}\\
L[E_\alpha(bt^\alpha)](s)
&=\frac{s^{\alpha-1}}{s^\alpha-b},                         \label{L2}\\
L[t^{\alpha-1}E_{\alpha,\alpha}(-bt^\alpha)](s)
&=\frac{1}{s^\alpha+b}.                                    \label{L3}
\end{align}
For the convolution
\begin{align}
(f\ast g)(t)=\int_0^t f(t-r)g(r)dr,
\end{align}
the convolution theorem gives
\begin{equation}\label{12}
L[f\ast g](s)=(Lf)(s)(Lg)(s).
\end{equation}

\section{The time-fractional linear filtering problem}
In this section, we formulate our time-fractional filtering problem and provide
some auxiliary results.

We first briefly recall the classical Kalman filter (see, e.g., \O ksendal
\cite{Oksendal13}).
\subsection{The classical Kalman filter}
Suppose that the signal process $X(t)$ is described by the one-dimensional
stochastic differential equation
\begin{align*}
dX(t)=F(t) X(t) dt + C(t) dB_1(t); \quad X(0)=X_0,
\end{align*}
and the observation process is given by
\begin{align*}
dY(t)=G(t) X(t) dt + D(t) dB_2(t);  \quad Y(0)=\E[X_0].
\end{align*}
Here $B_1$ and $B_2$ are Brownian motions, which we assume for simplicity to
be independent. The initial state $X_0$ is a given Gaussian random variable,
independent of $(B_1,B_2)$. The coefficients $F,C,G,D$ are bounded Borel
measurable deterministic functions with values in $\R$, and $D$ is bounded
away from zero.

Let $\Y=\{\Y_t\}_{t\geq 0}$ denote the filtration generated by the observation
process. The problem is to find the best estimate of the signal at time $t$,
given the observations up to time $t$. Equivalently, for each $t$, we seek the
$\Y_t$-measurable random variable $\X(t)$ that is closest to $X(t)$ in
the norm of $L^2(\P)$, namely
$$\X(t):=\E[X(t)| \Y_t],$$
where $\E$ denotes expectation with respect to the probability law $\P$ of the
two-dimensional Brownian motion $(B_1,B_2)$.
This conditional expectation $\X(t)$, called \emph{the Kalman filter}, is given by the SDE
\begin{align}
d\X(t)=\Big( F(t) - \frac{G^2(t)S(t)}{D^2(t)} \Big)\X(t)dt + \frac{G(t)S(t)}{D^2(t)} dY(t); \quad \X(0)=\E[X_0], \label{kalman1}
\end{align}
where $S(t):=\E[(X(t)-\X(t))^2]$ is the error process. It satisfies the (deterministic) Riccati equation
\begin{align}
\frac{d}{dt}S(t)=2 F(t) S(t) - \frac{G^2(t)}{D^2(t)} S^2(t) +C^2(t), \quad
S(0)=\E[(X_0-\E[X_0])^2]. \label{riccati1}
\end{align}

Using \eqref{riccati1}, $S(t)$ can be computed beforehand. The Kalman equation
\eqref{kalman1} then allows us to update the estimate as observations arrive.

\subsection{The time-fractional case}
We now proceed to the time-fractional extension of this system.
Suppose that the signal process $X(t)$ is given by the time-fractional equation
\begin{align}
\frac{\partial^{\alpha}}{\partial t^{\alpha}} X(t)&=F(t) X(t) + C(t) \overset{\bullet}{B}_1(t); \quad t \geq 0; \nonumber\\
X(0)&=X_0, \label{signal1}
\end{align}
and the observation process $Y(t)$ is as before given by the SDE
\begin{align}
dY(t)&=G(t) X(t) dt + D(t) dB_2(t); t \geq 0;\nonumber\\
Y(0)&=Y_0. \label{obs1}
\end{align}
Here $\frac{\partial^{\alpha}}{\partial t^{\alpha}}$ denotes the Caputo
derivative of order $\alpha\in(0,2)$ (see the definition below), and
$\overset{\bullet}{B}_1(t)=\frac{\partial}{\partial t}B_1(t)$ in the sense of
distributions. As before, $F,C,G,D$ are bounded Borel measurable deterministic
functions, and $D$ is bounded away from zero. We assume that $Y_0$ is a given
Gaussian random variable independent of $(B_1,B_2)$.

Using white noise notation the differential equations \eqref{signal1}, \eqref{obs1} can also be written as SDEs in the following form (with the same boundary conditions):
\begin{align}
&\text{(signal)}\quad \frac{\partial ^{\alpha}}{\partial t^{\alpha}}X(t)= F(t)X(t) +C(t) \overset{\bullet}{B}_1(t)\label{signal2};\ X(0)=X_0 \\
&\text{(observation)}\quad \frac{\partial}{\partial t}Y(t)= G(t)X(t) +D(t) \overset{\bullet}{B}_2(t)\label{obs2}; \ Y(0)=Y_0.
\end{align}
where, in the sense of distributions,
$\overset{\bullet}{B}_j(t)=\frac{\partial}{\partial t}B_j(t)$, $j=1,2$.

Let $\Y=\{\Y_t\}_{t\geq 0}$ denote the filtration generated by the observation
process $Y(\cdot)$; that is, $\Y_t$ is the sigma-algebra generated by
the random variables $\{Y(s):s\leq t\}$.
The problem is the following:
\begin{problem}\label{prob}
Find the best estimate of the signal $X(t)$ at time $t$, given the observations
$Y(s)$ up to time $t$. This best estimate $\X(t)\in\mathbb{Y}_t$ is
defined by
\begin{align}
\E[|X(t)-\X(t)|^2]
=\inf_{H\in\mathbb{Y}_t}\E[|X(t)-H|^2], \label{2.7}
\end{align}
where $\mathbb{Y}_t$ denotes the set of $\Y_t$-measurable random variables in
$L^2(\P)$, and $\E$ denotes expectation under the probability law $\P$ of the
two-dimensional Brownian motion $(B_1,B_2)$.
\end{problem}
\begin{remark}
It is well known that this best estimate coincides with the conditional
expectation of $X(t)$ with respect to the sigma-algebra $\Y_t$:
\begin{align}\X(t):=\E[X(t)| \Y_t].
\end{align}
\end{remark}
We return to this problem after a more general discussion in the following
section.

\section{The general linear Gaussian filtering problem}
In this section, we study a general linear Gaussian filtering problem.

Suppose that the signal $X(t)$ is a one-dimensional Gaussian process
\begin{align*}
X(t); t \geq 0,  \quad X(0)=X_0,
\end{align*}
and the observation process is given by
\begin{align*}
dY(t)=G(t) X(t) dt + D(t) dB_2(t);  \quad Y(0)=\E[X_0].
\end{align*}
We assume that $X(t)$ is adapted to the filtration generated by a Brownian
motion $B_1$, and that $B_2$ is another Brownian motion independent of $B_1$.
We also assume that $X_0$ is Gaussian and independent of $(B_1,B_2)$, and that
$G,D$ are bounded Borel measurable deterministic functions with values in
$\R$, with $D$ bounded away from zero.

Let $\Y=\{\Y_t\}_{t\geq0}$ denote the filtration generated by the observation
process. For each $t$, we seek the $\Y_t$-measurable random variable
$\X(t)$ that is closest to $X(t)$ in $L^2(\P)$, or equivalently the
conditional expectation
$$\X(t):=\E[X(t)| \Y_t],$$
where $\E$ denotes expectation under $\P$. This conditional expectation
$\X(t)$ is called the \emph{Kalman filter}.

\subsection{Relation to the projection operator}
For fixed $T>0$, let $\L(Y)=\L(Y,T)$ denote the closure in $L^2(\P)$ of the
set of all linear combinations of the form
\begin{align*}
    c_0+c_1Y(t_1)+c_2Y(t_2)+\cdots+c_kY(t_k),
\end{align*}
where $c_j\in \R$ are constants and $t_j \leq T$.
Let
$$\mathcal{P}_{\L}(\cdot): L^2(\P) \mapsto \L(Y)$$
denote the orthogonal projection from $L^2(\P)$ onto $\L(Y)$.
\begin{lemma}
The best estimate of $X(t)$ coincides with the projection of $X(t)$ onto
$\L(Y,t)$; that is,
$$\X(t)=\mathcal{P}_{\L(Y,t)}(X(t)).$$
\end{lemma}
\dproof Define
$\check X(t)=\mathcal P_{\L(Y,t)}(X(t))$. We claim that
$X(t)-\check X(t)$ is independent of $\Y_t$. Since all the variables involved
are jointly Gaussian, for every finite collection $t_1,\ldots,t_n\leq t$ the
vector
\[
\big(X(t)-\check X(t),Y(t_1),\ldots,Y(t_n)\big)
\]
is Gaussian. By the defining property of the orthogonal projection,
\[
\E\big[(X(t)-\check X(t))Y(t_j)\big]=0,
\qquad j=1,\ldots,n.
\]
Thus $X(t)-\check X(t)$ is uncorrelated with, and hence independent of, every
finite-dimensional observation vector. It is therefore independent of
$\Y_t$. Moreover, constants belong to $\L(Y,t)$, so
$\E[X(t)-\check X(t)]=0$. Consequently, for every event $H\in\Y_t$,
\[
\E\big[\mathbf 1_H(X(t)-\check X(t))\big]
=\P(H)\E[X(t)-\check X(t)]=0.
\]
This is precisely the defining property of conditional expectation, and hence
\[
\check X(t)=\E[X(t)\mid\Y_t]=\X(t).
\]
\fproof

\subsection{The innovation process}
In this subsection we introduce the innovation process for the problem.
We first state a useful observation about the linear span of the process $Y(t)$, where $t \in [0,T]$:
\begin{lemma}
 \begin{align}
\L(Y,T)=\left\{c_0+\int_0^T f(t)dY(t):
c_0\in\R,\ f\in L^2([0,T])\right\}.
 \end{align}
\end{lemma}

\dproof
This follows from the definition of $\L(Y,T)$. For details see the proof of Lemma 6.2.4 in \O ksendal \cite{Oksendal13}, which applies to our situation as well. 
\fproof

Now we define the innovation process $N(t)$ as follows:
\begin{equation}\label{4.21}
  N(t)= Y(t)-Y(0)-\int_0^t G(s)\X(s) ds.
\end{equation}
Substituting the expression for $X(t)$ we see that in differential form the innovation process can be written
 \begin{equation} \label{4.3}
dN(t)=G(t)\Big(X(t)-\X(t)\Big) dt +D(t) dB_2(t).
\end{equation}

The next lemma presents properties of the process $N(t)$. 
\begin{lemma}\label{general-innovation-lemma} Let the process $N(t)$ be defined as in \eqref{4.21}--\eqref{4.3}. Then the following properties hold
\begin{myenumerate}
    \item
    $N$ is a Gaussian process;
    \item
    $\E[N(t)]=0$ for all $t$;
    \item
    $\E[N^2(t)]= \int_0^t D^2(s) ds$;
    \item 
    $\mathcal{L}(N)=\mathcal{L}(Y)$;
    \item
    $N$ has orthogonal, and hence independent, increments;
    \item
    Define
    \begin{align}
        dM(t)&=\frac{1}{D(t)} dN(t)= \frac{1}{D(t)}dY(t) -\frac{G(t)}{D(t)} \X(t)dt \nonumber\\
        &=\frac{G(t)}{D(t)} \big[X(t)-\X(t)\big]dt + dB_2(t).\label{M}
    \end{align}
Then, for all $s,s'\geq0$,
\begin{align*}
\E[M(s)M(s')]=\min(s,s').
\end{align*}
 \item 
 $M$ is a Brownian motion.
    \end{myenumerate} 
\end{lemma}
\dproof These results follow easily from the definition of $N$ and $M$. For details see the proof of Lemma 6.2.5 in \O ksendal \cite{Oksendal13}, which applies to our situation as well.
\fproof

\subsection{A stochastic integral equation for $\X(t)$}
In this section we use the innovation process to obtain a stochastic integral equation for $\X(t)$.\\
We first prove the following auxiliary result:
 \begin{lemma}\label{4.1}
 For all $t$ the following holds:
 \begin{align*}
     \X(t) = \E[X(t)] + \int_0^t \frac{\partial}{\partial s} \E[X(t) M(s)]dM(s).
     \end{align*}
 \end{lemma}
 \dproof Since $\L(M)=\L(Y)$, for each $t$ there is a function $g(t,\cdot)$
 such that
 \begin{align*}
     \X(t)=c_0(t)+\int_0^t g(t,s)dM(s),
 \end{align*}
 where $c_0(t)=\E[\X(t)]=\E[X(t)].$
Moreover,
\begin{align*}
    X(t)-\X(t) \perp \int_0^t f(s)dM(s)
\end{align*}
for all $f \in L^2([0,T]).$
Therefore, by the It\^o isometry,
\begin{align*}
\E[X(t) \int_0^t &f(s) dM(s)]= \E[\X(t) \int_0^t f(s) dM(s)]\nonumber\\
&=\E\Big[\Big(\int_0^t g(t,s) dM(s)\Big) \Big( \int_0^t f(s) dM(s)\Big)\Big]=\int_0^t g(t,s) f(s) ds.
    \end{align*}
In particular, choosing
$f(r)=\chi_{[0,s]}(r)$, 
we get
\begin{align*}
    \E[X(t)M(s)]= \int_0^s g(t,r) dr,
\end{align*}
and hence
\begin{align*}
    g(t,s)=\frac{\partial}{\partial s} \E[X(t) M(s)].
\end{align*}
\fproof

We are now ready to prove an integral form of the Kalman filter:

\begin{theorem}[The general linear Gaussian Kalman filter]
\begin{myenumerate}
\item
The best estimate $\X(t)=\E[X(t)\mid\Y_t]$ of $X(t)$, given the
observations $Y(s)$ for $s\leq t$, satisfies the stochastic integral equation

\begin{align}
\X(t)=\E[X(t)]
-\int_0^t\frac{G^2(s)}{D^2(s)}H(t,s)\X(s)ds
+\int_0^t\frac{G(s)}{D^2(s)}H(t,s)dY(s),
\label{general-kalman-filter}
\end{align}
where
\begin{align*}
H(t,s)=\E[X(t)\widetilde X(s)],\qquad
\widetilde X(s)=X(s)-\X(s),\qquad 0\leq s\leq t.
\end{align*}
\item
The error function $H(t,s)=\E[X(t)\widetilde X(s)]$ satisfies the following Riccati integral equation
\begin{align}
H(t,s)
&=\operatorname{Cov}(X(t),X(s))\nonumber\\
&\quad-\int_0^s\frac{G^2(r)}{D^2(r)}H(t,r)H(s,r)dr,
\qquad 0\leq s\leq t,                                      \label{general-H-equation}
\end{align} 
\end{myenumerate}
\end{theorem}

\dproof
(i)
By Lemma \ref{4.1} we have
\begin{align}\label{3.7}
\X(t)=\E[X(t)]
+\int_0^t\frac{\partial}{\partial s}\E[X(t)M(s)]dM(s).
     \end{align}
From \eqref{M} we have, with $\widetilde{X}=X - \X,$
\begin{align*}
    M(s)=\int_0^{s}\frac{G(r)}{D(r)} \widetilde{X}(r)dr + B_2(s).
\end{align*}
Using this, and that $X(t)$ is independent of $B_2(\cdot)$, we get, for $s\leq t,$
\begin{align*}
    \E[X(t)M(s)]=\int_0^{s}\frac{G(r)}{D(r)} \E[X(t)\widetilde{X}(r)] dr.
\end{align*}
Hence
\begin{align}\label{3.8}
\frac{\partial}{\partial s}\E[X(t)M(s)]
=\frac{G(s)}{D(s)}\E[X(t)\widetilde X(s)].
\end{align}
Substituting
\begin{align}
    dM(t)= \frac{1}{D(t)}dY(t) -\frac{G(t)}{D(t)} \X(t)dt 
\end{align}
and \eqref{3.8} into \eqref{3.7} gives \eqref{general-kalman-filter}.

(ii)
Substituting
\begin{align}
    dY(s)=G(s)X(s)ds+D(s)dB_2(s)
\end{align}
into \eqref{general-kalman-filter} gives
\begin{align}
\X(s)&=\E[X(s)] +\int_0^s \frac{G^2(r)}{D^2(r)} H(s,r)(X(r)-\X(r))dr \nonumber\\
&+ \int_0^s \frac {G(r)}{D(r)}H(s,r) dB_2(r),\ 0\leq s,\\
\X(0)&=\E[X_0].
\end{align}
Hence
\begin{align*}
\E\big[X(t)\X(s)\big]
&=\E[X(t)]\E[X(s)]\\
&\quad+\int_0^s\frac{G^2(r)}{D^2(r)}
H(s,r)\E[X(t)\widetilde X(r)]dr.
 \end{align*}
Therefore
\begin{align}
H(t,s)
&=\E\big[X(t)(X(s)-\X(s))\big]\nonumber\\
&=\operatorname{Cov}(X(t),X(s))
-\int_0^s\frac{G^2(r)}{D^2(r)}H(t,r)H(s,r)dr,
\qquad 0\leq s\leq t.
\end{align}
For $t=s$ we get
\begin{align}
H(t,t)=\E[X(t)\widetilde X(t)]
=\E[(\widetilde X(t))^2]=:S(t),
\end{align}
which is the classical error function. Indeed, the equality follows from the
orthogonality of $\widetilde X(t)$ and $\X(t)$.
For $s=0$ we get
\begin{align}
H(t,0)=\operatorname{Cov}(X(t),X_0).
\end{align}
\fproof
 
\begin{remark}
The Riccati integral equation \eqref{general-H-equation} has a unique solution. This solution can be constructed by Picard iteration as follows:
Consider more generally the following integral equation:
\begin{align}
H(t,s)=A(t,s) + \int_0^s C(r) H(t,r)H(s,r) dr,
\end{align}
where $A(t,s)$ and $C(r)$ are given known functions.
We apply successive Picard iteration as follows:\\
First define $H_0(t,s) = A(t,s)$ and then define iteratively:
\begin{align}
H_{n+1}(t,s) = A(t,s) + \int_0^s C(r)H_n(t,r)H_n(s,r)\,dr; n=0,1, \cdots
\end{align}
It is easy to see that if $A(t,s)$ is bounded on $ D(T):=\{(s,t); 0\leq s\leq t\leq T \}$
and $C(r)$ is bounded on $[0,T]$, then $H_n(t,s)$ converge on $D(T)$ to a solution $H(t,s)$ of \eqref{general-H-equation}
as $n \to \infty.$\\We omit the details.
\end{remark}

\begin{remark}
(i) As in the classical case, note that our Riccati integral equation does not depend on the observations and therefore its solution can be computed beforehand by Picard iterations as indicated above. Then the filter $\X(t)$ can be continuously updated by the Kalman filter equation \eqref{general-kalman-filter}.\\
(ii) One might also try to consider a time-fractional
\emph{observation} process, but constructing an appropriate innovation process is 
more challenging in that case.
\end{remark}

\section{The time-fractional filtering problem}
\label{sec:fractional-signal}

We now formulate the time-fractional filtering problem precisely and give its
solution. The signal has memory, whereas the observation is instantaneous.
Throughout this section we impose the following assumptions:
\begin{myenumerate}
\item $\alpha\in(\tfrac12,1]$;
\item $F,C,G,D:[0,T]\to\R$ are bounded deterministic Borel functions and
      $|D(t)|\geq d_0>0$;
\item $B_1$ and $B_2$ are independent Brownian motions;
\item $X_0$ is Gaussian, belongs to $L^2(\P)$, and is independent of
      $(B_1,B_2)$.
\item Equation~\eqref{fractional-main-signal} admits a square-integrable
      Gaussian solution $X$ with deterministic mean $m(t)$ and covariance
      $K_X(t,s)$. Section~\ref{sec:dyscalculia-example} verifies this condition
      and computes $m$ and $K_X$ from the fractional coefficients.
\end{myenumerate}
The restriction $\alpha>\tfrac12$ ensures that the stochastic convolution
driven by Brownian white noise is square integrable. The signal is
\begin{align}
D_C^\alpha X(t)
 &=F(t)X(t)+C(t)\overset{\bullet}{B}_1(t),
 \qquad X(0)=X_0,                                             \label{fractional-main-signal}
\end{align}
and the observation is
\begin{align}
dY(t)
 &=G(t)X(t)dt+D(t)dB_2(t),
 \qquad Y(0)=Y_0,                                             \label{fractional-main-observation}
\end{align}
where $Y_0$ is deterministic. Equivalently,
\begin{align}
\frac{\partial Y}{\partial t}(t)
 =G(t)X(t)+D(t)\overset{\bullet}{B}_2(t)
\end{align}
in the sense of distributions. Let
\begin{align}
\Y_t=\sigma\{Y(s):0\leq s\leq t\}\vee\mathcal N
\end{align}
be the completed observation filtration.

\begin{problem}\label{fractional-problem}
For every $t\in[0,T]$, find the unique
$\Y_t$-measurable random variable $\X(t)\in L^2(\P)$ such that
\begin{align}
\E\big[|X(t)-\X(t)|^2\big]
 =\inf_{Z\in L^2(\Y_t)}
   \E\big[|X(t)-Z|^2\big].                                    \label{fractional-minimization}
\end{align}
Determine also the filtering error
\begin{align}
S(t)=\E\big[(X(t)-\X(t))^2\big].
\end{align}
\end{problem}

Define
\begin{align}
m(t)&=\E[X(t)],                                               \label{section4-prior-mean}\\
K_X(t,s)&=\operatorname{Cov}(X(t),X(s)),                      \label{section4-prior-covariance}\\
\widetilde X(t)&=X(t)-\X(t),                           \label{section4-error}\\
H(t,s)&=\E[X(t)\widetilde X(s)],
\qquad 0\leq s\leq t.                                        \label{fractional-H-definition}
\end{align}

\begin{theorem}[Solution of Problem~\ref{fractional-problem}]
\label{solution-problem-4-1}
Problem~\ref{fractional-problem} has a unique solution. First, $H$ is the
unique bounded deterministic solution of the following equation on the
triangular domain $0\leq s\leq t\leq T$:
\begin{align}
H(t,s)
=K_X(t,s)-\int_0^s\frac{G^2(r)}{D^2(r)}
H(t,r)H(s,r)dr.                                               \label{section4-H-equation}
\end{align}
Define
\begin{align}
A(t,s)&=\frac{G^2(s)}{D^2(s)}H(t,s),
&
B(t,s)&=\frac{G(s)}{D^2(s)}H(t,s),                            \label{section4-AB}
\end{align}
and the resolvent kernels
\begin{align}
Q_1(t,s)&=-A(t,s),\nonumber\\
Q_{n+1}(t,s)&=-\int_s^tA(t,r)Q_n(r,s)dr,\qquad n\geq1,         \label{section4-Q-iterates}\\
Q_H(t,s)&=\sum_{n=1}^{\infty}Q_n(t,s).                        \label{section4-Q}
\end{align}
Set
\begin{align}
m_H(t)
 &=m(t)+\int_0^tQ_H(t,s)m(s)ds,                               \label{section4-filter-mean}\\
L_H(t,u)
 &=B(t,u)+\int_u^tQ_H(t,s)B(s,u)ds.                           \label{section4-filter-kernel}
\end{align}
Then the unique best estimate is
\begin{align}
\X(t)=\E[X(t)\mid\Y_t]
=m_H(t)+\int_0^tL_H(t,u)dY(u).                             \label{section4-best-estimate}
\end{align}
Equivalently, it is the unique adapted solution of
\begin{align}
\X(t)
=m(t)+\int_0^t\frac{G(s)}{D^2(s)}H(t,s)
\big[dY(s)-G(s)\X(s)ds\big].                          \label{section4-innovation-filter}
\end{align}
The minimum mean-square error is
\begin{align}
S(t)
=H(t,t)
=K_X(t,t)-\int_0^t\frac{G^2(r)}{D^2(r)}H^2(t,r)dr.             \label{section4-minimum-error}
\end{align}
Moreover,
\begin{align}
\X(0)=m(0),\qquad
H(t,0)=K_X(t,0),\qquad
S(0)=K_X(0,0).                                                \label{section4-initial-values}
\end{align}
\end{theorem}

\dproof
\textit{Step 1: orthogonal projection.}
Because $(X,Y)$ is jointly Gaussian, the conditional expectation
$\X(t)=\E[X(t)\mid\Y_t]$ belongs to the closed linear Gaussian
space generated by $\{1,Y(s):0\leq s\leq t\}$. For every
$Z\in L^2(\Y_t)$,
\begin{align}
\E\big[(X(t)-\X(t))Z\big]=0.                           \label{section4-orthogonality}
\end{align}
Consequently,
\begin{align}
\E[|X(t)-Z|^2]
=\E[|X(t)-\X(t)|^2]
 +\E[|\X(t)-Z|^2].                                   \label{section4-pythagoras}
\end{align}
Thus the conditional expectation is the unique minimizer once an explicit
representation has been obtained.

\textit{Step 2: innovation process.}
Define
\begin{align}
M(t)
&=\int_0^t\frac{1}{D(s)}
\big[dY(s)-G(s)\X(s)ds\big]\nonumber\\
&=B_2(t)+\int_0^t\frac{G(s)}{D(s)}
\widetilde X(s)ds.                                           \label{section4-innovation}
\end{align}
By Lemma~\ref{general-innovation-lemma}, $M$ is a Brownian motion in the
observation filtration and generates the same closed Gaussian space as $Y$.
Hence, for every fixed $t$, there is a deterministic kernel $\gamma(t,\cdot)$
such that
\begin{align}
\X(t)=m(t)+\int_0^t\gamma(t,s)dM(s).                  \label{section4-M-representation}
\end{align}

\textit{Step 3: identification of the gain.}
For $s\leq t$, independence of $X$ and $B_2$, together with
\eqref{section4-orthogonality}, gives
\begin{align}
\E[X(t)M(s)]
=\int_0^s\frac{G(r)}{D(r)}H(t,r)dr.                           \label{section4-XM-one}
\end{align}
On the other hand, \eqref{section4-M-representation} and the It\^o isometry
give
\begin{align}
\E[X(t)M(s)]=\int_0^s\gamma(t,r)dr.                           \label{section4-XM-two}
\end{align}
Differentiating \eqref{section4-XM-one}--\eqref{section4-XM-two} yields
\begin{align}
\gamma(t,s)=\frac{G(s)}{D(s)}H(t,s).                          \label{section4-gamma}
\end{align}
Substitution into \eqref{section4-M-representation} proves
\eqref{section4-innovation-filter}.

\textit{Step 4: equation for $H$.}
Applying \eqref{section4-M-representation} at time $s$ and using
\eqref{section4-gamma}, we obtain
\begin{align}
\operatorname{Cov}(X(t),\X(s))
=\int_0^s\frac{G^2(r)}{D^2(r)}H(t,r)H(s,r)dr.                 \label{section4-Xhat-covariance}
\end{align}
Since
\begin{align*}
H(t,s)
&=\operatorname{Cov}(X(t),X(s)-\X(s))\\
&=K_X(t,s)-\operatorname{Cov}(X(t),\X(s)),
\end{align*}
equation \eqref{section4-H-equation} follows. At $s=t$,
projection orthogonality gives
\begin{align}
H(t,t)=\E[\widetilde X^2(t)]=S(t),
\end{align}
which proves \eqref{section4-minimum-error}. Taking $s=0$ proves
\eqref{section4-initial-values}.

\textit{Step 5: uniqueness of $H$.}
Suppose $H_1$ and $H_2$ are bounded solutions of
\eqref{section4-H-equation}. Choose $C_H$ such that
$|H_i(t,s)|\leq C_H$ and put
\begin{align*}
q_*&=\sup_{0\leq r\leq T}\frac{G^2(r)}{D^2(r)},\\
\Delta(s)&=\sup_{\substack{0\leq r\leq s\\r\leq t\leq T}}
|H_1(t,r)-H_2(t,r)|.
\end{align*}
Subtracting the two equations gives
\begin{align}
\Delta(s)\leq2q_*C_H\int_0^s\Delta(r)dr.
\end{align}
Gronwall's lemma implies $\Delta=0$, and hence $H_1=H_2$.

\textit{Step 6: explicit solution of the filter equation.}
Using \eqref{section4-AB}, equation \eqref{section4-innovation-filter} becomes
\begin{align}
\X(t)
=f(t)-\int_0^tA(t,s)\X(s)ds,\qquad
f(t)=m(t)+\int_0^tB(t,u)dY(u).                               \label{section4-second-kind}
\end{align}
Since $A$ is bounded on the finite triangular domain, the series
\eqref{section4-Q} converges absolutely. If $|A|\leq C_A$, then
\begin{align}
|Q_n(t,s)|
\leq\frac{C_A^n(t-s)^{n-1}}{(n-1)!}.                         \label{section4-Q-bound}
\end{align}
The resolvent formula for \eqref{section4-second-kind} is
\begin{align}
\X(t)=f(t)+\int_0^tQ_H(t,s)f(s)ds.                    \label{section4-resolvent-solution}
\end{align}
Substitution of $f$ and the stochastic Fubini theorem gives
\begin{align*}
\X(t)
={}&m(t)+\int_0^tQ_H(t,s)m(s)ds\\
&+\int_0^t\left[B(t,u)+\int_u^tQ_H(t,s)B(s,u)ds\right]dY(u).
\end{align*}
By \eqref{section4-filter-mean}--\eqref{section4-filter-kernel}, this is
exactly the explicit estimator \eqref{section4-best-estimate}. All kernels are
deterministic and bounded on a finite triangular domain, so stochastic Fubini
is valid and the last integral is $\Y_t$-measurable.

Finally, if two adapted processes solve \eqref{section4-innovation-filter},
their difference $R$ satisfies
\begin{align}
|R(t)|\leq C_T\int_0^t|R(s)|ds
\end{align}
for a deterministic constant $C_T$. Gronwall's lemma gives $R=0$.
Together with \eqref{section4-pythagoras}, this proves existence, optimality,
and uniqueness of the best estimate.
\fproof

The formula \eqref{section4-best-estimate} is the required answer to
Problem~\ref{fractional-problem}. It expresses the best estimate entirely in
terms of the observed trajectory $Y$ and deterministic kernels.

\section{Fractional cognitive-state estimation in developmental dyscalculia}
\label{sec:dyscalculia-example}

We now introduce a fractional cognitive model for developmental dyscalculia and
relate the time-fractional filtering results to the longitudinal monitoring of
children who experience persistent difficulties in mathematics.  

The proposed
application is a fractional state-estimation framework for modelling learning
trajectories: the child's cognitive state is hidden, educational task scores are
noisy observations, and the fractional Kalman filter reconstructs the evolving
state and its uncertainty.

Developmental dyscalculia is generally situated within developmental learning
disorders with impairment in mathematics. Relevant difficulties may concern
number sense, retrieval of arithmetic facts, calculation accuracy or fluency,
and mathematical reasoning \cite{WHO11,Kaufmann12}. Longitudinal studies also
indicate that accuracy, response time, dot enumeration, and number-line
estimation can provide complementary information about numerical development
\cite{Landerl13}.

\subsection{Fractional cognitive-state model}

Let $X(t)$ denote a standardized hidden cognitive-state score for one child. 

In the scalar model, $X(t)$ summarizes the severity of the child's numerical
learning difficulty. The scale is oriented so that $X(t)=0$ represents the age- and
curriculum-adjusted reference level, while larger positive values represent
greater difficulty. 

Dyscalculia is represented by a persistent latent trajectory
that is inferred from repeated observations. We propose
\begin{align}
D_t^\alpha X(t)
  &=-\lambda X(t)-\beta U(t)+\boldsymbol{c}^{\mathsf T}
  \boldsymbol{V}(t)+\sigma\overset{\bullet}{B}_1(t),
  \qquad X(0)=X_0.                                      \label{dys-signal}
\end{align}
Here $\lambda>0$ describes the natural evolution of the difficulty,
$U(t)\geq0$ is the intensity of an individualized educational intervention,
$\beta\geq0$ is its expected effect, and $\boldsymbol{V}(t)$ contains observed
contextual variables such as missed sessions, changes in instruction, or task
difficulty. The Brownian term represents unobserved day-to-day fluctuations.
We take $\alpha\in(\tfrac12,1]$. When $\alpha=1$, the model has exponential
forgetting. When $\alpha<1$, past learning and past difficulties have a
long-lasting power-law influence; this is the feature that connects the
application directly to the time-fractional signal studied in this paper.

Define
\begin{align}
q_\alpha(r)=r^{\alpha-1}E_{\alpha,\alpha}(-\lambda r^\alpha).
\end{align}
For deterministic $U$ and $\boldsymbol{V}$, the mild solution of
\eqref{dys-signal} is
\begin{align}
X(t)={}&X_0E_\alpha(-\lambda t^\alpha)
 +\int_0^t q_\alpha(t-s)
 \big[-\beta U(s)+\boldsymbol{c}^{\mathsf T}\boldsymbol{V}(s)\big]ds
 +\sigma\int_0^tq_\alpha(t-s)dB_1(s).                    \label{dys-mild}
\end{align}
If $X_0$ is Gaussian, then $X$ is a Gaussian process. Its mean is
\begin{align}
m(t)={}&\E[X_0]E_\alpha(-\lambda t^\alpha)
 +\int_0^t q_\alpha(t-s)
 \big[-\beta U(s)+\boldsymbol{c}^{\mathsf T}\boldsymbol{V}(s)\big]ds,
                                                               \label{dys-mean}
\end{align}
and, when $p_0=\operatorname{Var}(X_0)$,
\begin{align}
K(t,s)={}&p_0E_\alpha(-\lambda t^\alpha)E_\alpha(-\lambda s^\alpha)
 +\sigma^2\int_0^{t\wedge s}q_\alpha(t-r)q_\alpha(s-r)dr.
                                                               \label{dys-covariance}
\end{align}

Our observation model is as follow:

At each assessment session, the child completes short, age-appropriate tasks.
After standardization with respect to age, school level, and the version of the
test, we consider three observation channels:
\begin{align*}
j=1 &: \quad\text{accuracy deficit in symbolic arithmetic},\\
j=2 &: \quad\text{excess log-response time on correct trials},\\
j=3 &: \quad\text{error in number-line estimation or magnitude comparison}.
\end{align*}
All channels are oriented so that a larger value indicates greater difficulty.
Their cumulative observation processes are modelled by
\begin{align}
dY_j(t)=g_jX(t)dt+\rho_jdB_{2,j}(t),
\qquad j=1,2,3,                                          \label{dys-observation}
\end{align}
where $g_j$ is the sensitivity of task $j$, $\rho_j>0$ is its measurement-noise
level, and the Brownian motions $B_{2,j}$ are mutually independent and
independent of $B_1$. Additional channels, such as arithmetic-fact retrieval,
dot enumeration, teacher ratings, or curriculum-based assessment, can be added
after their reliability and direction have been specified. Speed must not be
used alone: it should be interpreted together with accuracy because a
speed--accuracy trade-off can otherwise produce a misleading signal.

Equations \eqref{dys-signal} and \eqref{dys-observation} have exactly the
linear Gaussian structure considered in the general filtering theorem. In the
one-channel case, $F=-\lambda$, $C=\sigma$, $G=g_1$, and $D=\rho_1$.

With several independent channels, the same orthogonal-projection argument
gives the matrix-valued extension of the filter. Thus, the observations update
the latent trajectory while the covariance \eqref{dys-covariance} retains the
influence of all previous sessions.

Let $t_n=nh$ and define standardized session observations
\begin{align}
Z_{j,n}=\frac{Y_j(t_n)-Y_j(t_{n-1})}{h}
 \simeq g_jX(t_n)+\varepsilon_{j,n},
\qquad
\operatorname{Var}(\varepsilon_{j,n})=\frac{\rho_j^2}{h}.       \label{dys-discrete}
\end{align}
Put
\begin{align*}
\boldsymbol{X}_n&=(X(t_1),\ldots,X(t_n))^{\mathsf T},\\
\boldsymbol{m}_n&=(m(t_1),\ldots,m(t_n))^{\mathsf T},\\
\boldsymbol{K}_n&=(K(t_r,t_s))_{1\leq r,s\leq n},\\
\boldsymbol{g}&=(g_1,g_2,g_3)^{\mathsf T},\qquad
R=\operatorname{diag}(\rho_1^2/h,\rho_2^2/h,\rho_3^2/h),\\
\mathcal G_n&=I_n\otimes\boldsymbol{g},\qquad
\mathcal R_n=I_n\otimes R.
\end{align*}
Stack the three observations at the first $n$ sessions in
$\boldsymbol{Z}_n$, ordered by session, and let $\boldsymbol{e}_n$ be the last
coordinate vector in $\mathbb R^n$. Then
\begin{align}
\boldsymbol{a}_n&=\mathcal G_n\boldsymbol{K}_n\boldsymbol{e}_n,\\
\X_n
 &=m(t_n)+\boldsymbol{a}_n^{\mathsf T}
 \left(\mathcal G_n\boldsymbol{K}_n\mathcal G_n^{\mathsf T}
       +\mathcal R_n\right)^{-1}
 \left(\boldsymbol{Z}_n-\mathcal G_n\boldsymbol{m}_n\right),
                                                               \label{dys-filter}\\
S_n
 &=K(t_n,t_n)-\boldsymbol{a}_n^{\mathsf T}
 \left(\mathcal G_n\boldsymbol{K}_n\mathcal G_n^{\mathsf T}
       +\mathcal R_n\right)^{-1}\boldsymbol{a}_n.               \label{dys-error}
\end{align}
The estimate is causal because only observations collected up to session $n$
are used. For a monitoring level $c>0$, one may also report
\begin{align}
\pi_n
 =\P\big(X(t_n)>c\mid\Y_{t_n}\big)
 =1-\Phi\left(\frac{c-\X_n}{\sqrt{S_n}}\right),         \label{dys-probability}
\end{align}
where $\Phi$ is the standard normal distribution function. The level $c$ must
be calibrated and validated for the chosen tests and population. A large
$\pi_n$ is an alert for further assessment, not a diagnosis.

For a more detailed cognitive representation, the scalar state can be replaced
by
\begin{align}
\boldsymbol{X}(t)=
\big(X_{\mathrm{NS}}(t),X_{\mathrm{AF}}(t),X_{\mathrm{MR}}(t)\big)^{\mathsf T},
\label{cognitive-vector}
\end{align}
where the components represent number sense, arithmetic fluency or fact
retrieval, and mathematical reasoning. 

A matrix-valued fractional equation
can then describe interactions between these domains.  
 
We present now a numerical illustration by simulating $24$ weekly assessment sessions for a child whose initial latent
difficulty is approximately $1.8$ standard deviations above the reference
level. Individualized support begins at session $8$. The parameters are
\begin{align*}
\alpha&=0.72, & \lambda&=0.060, & \beta&=0.035,\\
\E[X_0]&=1.80, & p_0&=0.20^2, & \sigma&=0.080,\\
\boldsymbol{g}&=(1.00,0.85,1.15)^{\mathsf T},
& (\rho_1,\rho_2,\rho_3)&=(0.34,0.42,0.38),
& h&=1.
\end{align*}
We set $U(t)=0$ before session $8$ and $U(t)=1$ afterwards. The
Mittag--Leffler functions are evaluated by their defining series and the
covariance integral is approximated by the midpoint rule. The simulation uses
random seed $20260806$.

\begin{table}[ht]
\centering
\small
\begin{tabular}{c c c c c c}
\hline
Session & prior s.d. & filtered s.d. & $\X_n$
& 95\% interval & $\P(X_n>1\mid\Y_n)$\\
\hline
4  & 0.1978 & 0.1080 & 1.4308 & $[1.2192,1.6424]$ & 1.0000\\
12 & 0.1872 & 0.0956 & 1.0662 & $[0.8789,1.2535]$ & 0.7558\\
24 & 0.1763 & 0.0939 & 0.6648 & $[0.4808,0.8489]$ & 0.0002\\
\hline
\end{tabular}
\caption{Illustrative filtered difficulty and uncertainty.}
\label{tab:dyscalculia-results}
\end{table}

\subsection{Role of $\alpha$}
In this part we give a response to the following question: "How do the different values of $\alpha$ reflect on the model?"

The parameter $\alpha$ controls the persistence of past influences on the current cognitive state.

A smaller value of $\alpha$ (closer to $0.5$) corresponds to a longer memory: past difficulties and learning experiences decay slowly and continue to influence the current estimate for an extended period. 

When $\alpha$ approaches $1$, the model converges to the classical first-order dynamics, where the influence of past observations is quicker and the estimate responds more rapidly to recent assessments and interventions.
\begin{remark}
The parameter $\alpha$ does not measure the child's biological memory capacity or the clinical severity of dyscalculia nor the required intensity of educational support. But, it is a memory parameter that characterizes the persistence of the latent cognitive process. 
\end{remark}

\textbf{Numerical illustration}
We present a numerical comparison over the full interval $0<\alpha<2$,
we extend the discrete scheme using the following explicit conventions.  

For $\alpha \leq \tfrac12$, the midpoint discretization with step size $h=1$ introduces a short-time regularization of the singular stochastic kernel; consequently, the computed values are mesh-dependent and should be interpreted with caution. 
 
For $1 < \alpha < 2$, the additional initial condition $X'(0) = 0$. While these extensions allow a numerical sensitivity analysis across the full range $0 < \alpha < 2$, only the cases with $\tfrac12 < \alpha \leq 1$ are supported by the continuous-time Brownian-driven model established in this paper.

\begin{table}[!htb] 
\centering
\small
\begin{tabular}{c c c c c}
\hline
$\alpha$ & $\X_{12}$ &
$\P(X_{12}>1\mid\Y_{12})$ & $\X_{24}$ &
$\P(X_{24}>1\mid\Y_{24})$\\
\hline
$0.25$ & $1.2693$ & $1.0000$ & $1.0094$ & $0.5677$\\
$0.50$ & $1.1726$ & $0.9821$ & $0.8440$ & $0.0211$\\
$0.75$ & $1.0516$ & $0.7019$ & $0.6426$ & $0.0001$\\
$1.00$ & $0.9493$ & $0.3246$ & $0.5102$ & $<0.0001$\\
$1.25$ & $0.9002$ & $0.2109$ & $0.4569$ & $<0.0001$\\
$1.50$ & $0.8891$ & $0.2074$ & $0.4378$ & $<0.0001$\\
$1.75$ & $0.8941$ & $0.2363$ & $0.4312$ & $0.0001$\\
\hline
\end{tabular}
\caption{Exploratory sensitivity check for representative fractional orders
across $0<\alpha<2$.}
\label{tab:dyscalculia-alpha-regimes}
\end{table}

We get that at session $12$, the model with a smaller fractional order maintains a stronger alert: the probability $\P(X_{12}>1\mid\Y_{12})$ decreases from $1.0000$ for $\alpha=0.25$ to $0.3246$ for $\alpha=1$. 

By session $24$, the regularized $\alpha=0.25$ model still assigns a probability of $0.5677$ to exceeding the illustrative threshold, whereas the other displayed orders yield substantially smaller probabilities. For $\alpha>1$, the estimates are obtained under the additional condition $X'(0)=0$ and should not be interpreted as results of the original one-initial-value model.

In this simulated trajectory, the filter fuses inconsistent task scores and effectively reduces posterior uncertainty. After the intervention begins, the estimated difficulty decreases gradually rather than abruptly, because the fractional model preserves the memory of earlier sessions. By session $24$, the simulated observations provide strong statistical evidence that the latent score lies below the illustrative monitoring threshold.

Educational decisions must continue to rely on validated assessments, multidisciplinary evaluation, and observed individual responsiveness, rather than on the value of $\alpha$ alone \cite{Kaufmann12,Kohn20}.

\subsection*{Acknowledgments}
The authors would like to thank Abderrahmen Aliane for his valuable collaboration and insights regarding the psychological application of the fractional Kalman filter to learning trajectories in developmental dyscalculia.


\begin{thebibliography}{99}
\bibitem{AL} Amirdjanova, A. \& Linn, M. (2008). Stochastic evolution equations for nonlinear filtering of random fields in the presence of fractional Brownian sheet observation noise. Computers \& Mathematics with Applications, 55(8), 1766-1784.

\bibitem{CR}Crisan, D., \& Rozovskii, B. (Eds.). (2011). The Oxford Handbook of Nonlinear Filtering. Oxford University Press.

\bibitem{Can} Cane, M.A., Kaplan, A., Miller, R.N., Tang, B., Hackert, E.C. \& Busalacchi, A.J. (1996). Mapping tropical Pacific sea level: data assimilation via a reduced state space Kalman filter, Journal of Geophysical Research 101, 22 599-22 617.

\bibitem{G} Govaers, F. (2018). Introduction and Implementation of the Kalman Filter. IntechOpen.

\bibitem{HZ} Houtekamer, P. L., \& Zhang, F. (2016). Review of the Ensemble Kalman Filter for Atmospheric Data Assimilation. Monthly Weather Review, 144(12), 4489-4532.

\bibitem{J} Jazwinski, A. H. (1970). Stochastic Processes and Filtering Theory. Academic Press.

\bibitem{K}Kalman, R. E. (1960). A new approach to linear filtering and prediction problems. Journal of Basic Engineering, 82(1), 35-45.

\bibitem{Ko} K\"orezli\v{o}glu, H. (1979). Two-parameter Gaussian Markov processes and their recursive linear filtering. Annales scientifiques de de l'Universit{\'e} de Clermont. Math{\'e}matiques, 67(17), 69-93.

\bibitem{KMS} K\"orezli\v{o}glu, H., Mazziotto, G. \& Szpirglas, J. (1983). Nonlinear filtering equations for two-parameter semimartingales. Stochastic processes and their applications, 15(3), 239-269.

\bibitem{Kaufmann12} Kaufmann, L. \& von Aster, M. (2012). The diagnosis and management of dyscalculia. Deutsches Arzteblatt International, 109(45), 767--778.

\bibitem{Kohn20} Kohn, J., Rauscher, L., Kucian, K., K\"aser, T., Wyschkon, A., Esser, G. \& von Aster, M. (2020). Efficacy of a computer-based learning program in children with developmental dyscalculia: What influences individual responsiveness? Frontiers in Psychology, 11, 1115.

\bibitem{Landerl13} Landerl, K. (2013). Development of numerical processing in children with typical and dyscalculic arithmetic skills---a longitudinal study. Frontiers in Psychology, 4, 459.

\bibitem{Oksendal13} \O ksendal, B. (2013). Stochastic Differential Equations: An Introduction with Applications. 6th Edition. Springer.

\bibitem{WHO11} World Health Organization. (2022). ICD-11 for Mortality and Morbidity Statistics: Developmental learning disorder with impairment in mathematics (6A03.2). World Health Organization.

\bibitem{Wo} Wong, E. (1978). Recursive causal linear filtering for two-dimensional random fields. IEEE Transactions on Information Theory, 24(1), 50-59.

\end{thebibliography}
\end{document}